\documentclass[a4paper,11pt]{amsart}
\usepackage{t1enc,amssymb,epsfig}
\usepackage[english]{babel}

\input xy
\xyoption{all}

\newtheorem{proposition}{Proposition}[section]
\newtheorem{theorem}[proposition]{Theorem}
\newtheorem{lemma}[proposition]{Lemma}
\newtheorem{definition}[proposition]{Definition}

\newtheorem{remark}[proposition]{Remark}

\begin{document}
\title[Iterating the hessian]{\textbf{Iterating the hessian: a dynamical
    system\\ on the moduli space of elliptic curves \\ and dessins d'enfants}}
\author{\sc Patrick Popescu-Pampu}
\address{Univ. Paris Diderot - Paris 7, Inst. de
  Maths.-UMR CNRS 7586, {\'e}quipe "G{\'e}om{\'e}trie et dynamique" \\case
  7012,  75205 Paris Cedex 13, France.}
\email{ppopescu@math.jussieu.fr}

\date{30.01.2009}

\thispagestyle{empty}
\begin{abstract}
  Each elliptic curve can be embedded uniquely in the projective
  plane, up to projective equivalence. The hessian curve of the
  embedding is generically a new elliptic curve, whose isomorphism
  type depends only on that of the initial elliptic curve. One gets like this
  a rational map from the moduli space of elliptic curves to itself. We call it
  \emph{the hessian dynamical system}. We compute it in terms of the
  $j$-invariant of elliptic curves. We deduce that, seen as a map from
  a projective line to itself, it has 3 critical values, which
  correspond to the point at infinity of the moduli space and to the
  two elliptic curves with special symmetries. Moreover, it sends the
  set of critical values into itself, which shows that all its
  iterates have the same set of critical values. One gets like this a
  sequence of dessins d'enfants. We describe an algorithm allowing to
  construct this sequence.  
\end{abstract}

\maketitle

\par\medskip\centerline{\rule{2cm}{0.2mm}}\medskip
\setcounter{section}{0}

\section{Introduction}

Consider a complex projective plane $\mathbb{P}^2$ with homogeneous 
coordinates \linebreak $[x:y:z]$. 

To each homogeneous polynomial $f \in \mathbb{C}[x,y,z]$ one can
associate the corresponding \emph{hessian polynomial} $\mathrm{Hess}(f)$, 
defined as the determinant:
\begin{equation} \label{hessian}
   \mathrm{Hess}(f) := \left | \begin{array}{ccc}
                         f_{xx} & f_{xy} & f_{xz}\\
                         f_{yx} & f_{yy} & f_{yz}\\
                         f_{zx} & f_{zy} & f_{zz}
                      \end{array} \right|.
\end{equation}
 
If $C:=Z(f)$ denotes the projective curve defined by the homogeneous
polynomial $f$ and $\mathrm{Hess}(C)$ denotes the curve defined by
$\mathrm{Hess}(f)$,   
one knows that the intersection points of $C$ and $\mathrm{Hess}(C)$ are
exactly the inflection points of $C$, which shows that this set has an
invariant meaning: it depends only on the couple $(\mathbb{P}^2, C)$,
and not on the chosen defining homogeneous polynomial. But more is
true: \emph{the 
whole curve $\mathrm{Hess}(C)$ is invariantly attached to $(\mathbb{P}^2,
C)$}. This remark was the starting point of the present study. 

If $\deg(f) =n \geq 1$, then $\deg(\mathrm{Hess}(f))=3(n-2)$. This shows that 
$C$ and $\mathrm{Hess}(C)$ have the same degree only when $n=3$. Restricting to
this case, we get a map from the space of projective plane cubic
curves to itself, which is equivariant with 
respect to the action of the group of projectivities of
$\mathbb{P}^2$. This shows that the map $C \rightarrow \mathrm{Hess}(C)$
descends to a rational map $H$ from the 
quotient of the space of smooth cubic curves by the group of
projectivities to itself. But this quotient space is 
\emph{the coarse moduli 
  space} $\mathcal{M}_1$ of elliptic curves (see Edidin \cite{E
  03}). The previous map extends to the compactification
$\overline{\mathcal{M}}_1 \simeq \mathbb{P}^1$. 
I propose: 

\begin{definition} \label{defdyn}
    The algebraic map $H: \overline{\mathcal{M}}_1 \rightarrow
    \overline{\mathcal{M}}_1$ which 
    associates to each elliptic curve the isomorphism type of the
    hessian of a smooth plane cubic curve corresponding to it is called 
    \textbf{the hessian dynamical system}.
\end{definition}

The aim of this paper is to compute the hessian dynamical system and
to start its dynamical study. I
believe that considering it could bring new insights into the theory
of elliptic curves. One could examine for example the relation between
dynamically defined subsets of $\mathcal{M}_1$  and the various
subsets with arithmetical meaning. More importantly, I believe 
that similar 
considerations related to higher dimensional classical invariant theory would
allow to construct higher dimensional dynamical systems with special
properties. 

Let me describe briefly the content of the paper. 
In Section \ref{Secinvj} are recalled various normal forms for plane
cubic curves and for each one of them, the expression of the classical
$j$-invariant. In Section \ref{Secomp} are computed the expression of
the hessian dynamical system $H$ in terms of the $j$-invariant (see
Theorem \ref{exphess}). In
Section \ref{Secdes} is showed that all the iterates of $H$ have $3$
critical values (see Proposition \ref{iterates}), which allows to
introduce an associated sequence 
of dessins d'enfants $(\Gamma_n)_{n \geq 1}$. Finally, in Section
\ref{Secomb}, is given an 
algorithm which allows to construct up to topological conjugacy the
sequence of preimages of the real axis by the iterates of $H$. In
particular, one gets a sequence of graphs in which the sequence of
dessins d'enfants introduced before embeds canonically, which gives us
an algorithm for constructing the sequence $(\Gamma_n)_{n \geq 1}$
(see Proposition \ref{embgraph}).

\medskip
While this paper was refereed, I learned that Pilgrim \cite{P 00} had
studied in general the relation between complex dynamics in dimension
one and dessins d'enfants and that Artebani \& Dolgachev have surveyed
in \cite{AD 06} the classical geometry of the Hesse pencil of
cubics. In Remark 3.5 of this last paper, they notice that the map
$\tilde{H}$ (see formula (\ref{hescub})) is an ``interesting example
of complex dynamics in one complex variable'', and that it was studied
from this view-point by Hollcroft \cite{H 26}. See also Remark
\ref{add} for other relations with the litterature.

\medskip

\section{The $j$-invariant of an elliptic curve} \label{Secinvj}

The classical $j$-invariant of an elliptic curve may be defined as
follows (see Hartshorne \cite[IV.4]{H 77}):

\begin{definition} \label{invj}
  Denote by $E_{\lambda}$ the elliptic curve birationally equivalent
  to the plane affine cubic with equation :
  $$ y^2 = x(x-1)(x-\lambda).$$
  Then its \textbf{$j$-invariant} is defined by:
  $$j(E_{\lambda}):= 2^8 \cdot
  \dfrac{(\lambda^2-\lambda+1)^3}{\lambda^2(\lambda-1)^2}. $$
\end{definition} 

The previous expression is adapted to the computation of the
$j$-invariant of an elliptic curve seen as the double cover of a
projective line 
$\mathbb{P}^1$, ramified over $4$ distinct points, the cross-ratio of
those points being $\lambda$. Indeed, $E_\lambda$ is the total space of
that covering over the complement of one of the $4$ points of
$\mathbb{P}^1$. 

If an elliptic curve is presented in Weierstarss normal form, one has
the following expression for its $j$-invariant (see Hartshorne
\cite[page 327]{H 77}):

\begin{proposition} \label{weij}
   Denote by $E_{g_2, g_3}$ the elliptic curve birationally equivalent
   to the smooth plane affine curve with equation:
   $$y^2 = 4x^3 -g_2  x -g_3,$$
   where $g_2 ^3 - 27 g_3^2 \neq 0$. Then its $j$-invariant is given
   by:
    $$j(E_{g_2, g_3})= 1728\cdot \dfrac{g_2^3}{g_2 ^3 - 27 g_3^2}.$$
\end{proposition}

In the sequel, we will work rather with the following normal form, for
reasons explained in the next section: 
\begin{equation} \label{linsyst}
   X_0^3 +X_1^3 +X_2^3 -3m \: X_0 X_1 X_2 =0.
\end{equation}

Denote by $C_m$ the plane projective cubic curve defined by equation
(\ref{linsyst}). As it is not easy to find a reference for the
following proposition, we add an elementary proof.

\begin{proposition} \label{projj}
  The $j$-invariant of the elliptic curve $C_m$ is given by:
  $$ j(C_m)= 27 \cdot \bigg(\dfrac{m(m^3+8)}{m^3 -1}\bigg)^3.$$
\end{proposition}

\textbf{Proof:} Define :
\begin{equation} \label{Jdef}
  J(C_m):= \bigg(\dfrac{m(m^3+8)}{m^3 -1}\bigg)^3.
\end{equation}

By Brieskorn \& Kn{\"o}rrer \cite[page 302]{BK 86}, we know that $J(C_m)$
and $j(C_m)$ are proportional, that is, there
exists $t \in \mathbb{C}$ such that:
\begin{equation} \label{relJj}
  J(C_m)= t \cdot j(C_m).
\end{equation}

In order to find $t$, it is enough to specialize (\ref{relJj}) to a 
cubic curve $C_m$ for which one knows how to compute both $J$ and
$j$. This is possible if one knows how to write $C_m$ in Weierstrass
normal form by a coordinate change. 

As the projectivisation of an
affine cubic in Weierstrass normal form has the property that the line
at infinity is tangent to it at an inflection point, we naturally begin by
choosing as line at infinity a tangent to $C_m$ at an inflection point. 

We choose the inflection point $(1 : -1 :0) \in C_m$. The tangent to
$C_m$ at this point has the equation:
  $$ X_0 + X_1 + mX_2=0.$$
We make then the following change of projective coordinates:
$$\left\{ \begin{array}{l}
           Z= X_0 + X_1 + m X_2\\
           U= \dfrac{1}{2}(X_0 + X_1)\\
           V= \dfrac{1}{2}(X_0 - X_1)
         \end{array} \right. \Longleftrightarrow 
  \left\{ \begin{array}{l}
            X_0 = U+V\\
            X_1 = U-V\\
            X_2 =\dfrac{1}{m}(Z-2U) 
          \end{array}  \right.  .$$
The equation (\ref{linsyst}) is transformed in:
  $$8 U^3 + \dfrac{1}{m^3}(Z-2U)^3 - 3(U^2 -V^2)Z=0.$$
By passing to the affine coordinates $u =U/Z, \: v= V/Z$, we find the
affine equation:
  $$ v^2 = \dfrac{8}{3}(\dfrac{1}{m^3}-1)u^3 + 
         (1-\dfrac{4}{m^3})u^2 + \dfrac{2}{m^3}u - \dfrac{1}{3m^3}.$$
We specialize now to $m = \sqrt[3]{4}$. After the new change of
variables $u= -2u_1, \: v= 2v_1$, the previous equation becomes:
 $$ v_1^2 =4 u_1 ^3 -\dfrac{1}{4}u_1 -\dfrac{1}{48}.$$

This shows that:
  $$ C_{\sqrt[3]{4}}\simeq E_{\frac{1}{4}, \frac{1}{48}}. $$
Combining this with relation (\ref{relJj}), we get:
  $$ J(C_{\sqrt[3]{4}}) = t \cdot j(E_{\frac{1}{4},\frac{1}{48}}).$$ 
From Proposition \ref{weij} and equation (\ref{Jdef}), we deduce that $t
=\dfrac{1}{27}$. \hfill $\Box$

\medskip

\section{Computation of the hessian dynamical 
  system}   \label{Secomp}

Consider again the 1-dimensional linear
system of plane projective cubics defined by the equation 
 (\ref{linsyst}).

The parameter $m \in \mathbb{C}$ is seen as an affine coordinate of the
projective line parametrizing the cubics of the pencil, the
homogeneous coordinates being $[1:m]$. 
To $m=\infty$ corresponds the cubic with equation
$X_0X_1X_2=0$, the union of the edges of the fundamental triangle in the
projective plane with the fixed homogeneous coordinates. 

An immediate computation shows that :
 \begin{equation}\label{hesscurve}
   \mathrm{Hess}(C_m)= C_{\frac{4-m^3}{3m^2}}.
 \end{equation}

\begin{remark}
  Equation (\ref{hesscurve}) shows that the hessian curve of a cubic
  expressed in the normal form (\ref{linsyst}) is again a cubic of the
  same normal form in the same system of homogeneous coordinates.
    If one had started instead from the expressions of Definition
    \ref{invj} or of Proposition \ref{weij}, one wouldn't have got
    expressions of the same normal form. This is the reason why we preferred
    to work with the normal form (\ref{linsyst}).
\end{remark}

One gets like this  a dynamical system $\tilde{H}: \mathbb{P}^1_{[1:m]} 
\rightarrow \mathbb{P}^1_{[1:m]}$, expressed by :
\begin{equation} \label{hescub}
      \tilde{H}(m) = \frac{4-m^3}{3m^2}
\end{equation}
where $\mathbb{P}^1_{[1:m]}:= \mathbb{A}^1_m \cup \{\infty\}$ denotes
the projective line obtained by adding one point at $\infty$ to the
affine line with coordinate $m$.

One has the following relation between the maps $\tilde{H}:
\mathbb{P}^1_{[1:m]}  \rightarrow \mathbb{P}^1_{[1:m]} ,\:  
H:\mathbb{P}^1_{[1:j]}  \rightarrow \mathbb{P}^1_{[1:j]},\:  
j:\mathbb{P}^1_{[1:m]}  \rightarrow \mathbb{P}^1_{[1:j]}$ :
$$ H \circ j= j\circ \tilde{H}.$$
Using Proposition \ref{projj}, it can be rewritten more explicitly as :
\begin{equation} \label{commdiag}
   H\bigg(27\bigg(\dfrac{m(m^3+8)}{m^3 -1}\bigg)^3\bigg) =
   27 \bigg(\dfrac{\tilde{H}(m)(\tilde{H}(m)^3+8)}{\tilde{H}(m)^3 -1}\bigg)^3 .
\end{equation}

By specializing the previous equality at $m\in \{ 0,1,
4^{\frac{1}{3}} \}$, we see that $H(0)=H(\infty)=\infty$ and
$H(2^8 \cdot 3^3)=0$. Moreover, when we see it as a holomorphic map from
$\mathbb{P}^1_{[1:m]}$ to itself, $\tilde{H}$ is of degree $3$, which
implies that $H$ is of degree $3$. As an immediate consequence of
these facts, we get:

\begin{lemma} \label{prelform}
  The rational fraction $H^*(j)$ of the variable $j$ is of the form:
$$ H^*(j)= \frac{(j-2^8\cdot 3^3)(j^2 + \alpha j +\beta)}{j(\gamma j +
  \delta)}$$
where $\alpha, \beta, \gamma, \delta \in \mathbb{C}$, $\beta \neq 0$
and $(\gamma, \delta) \neq (0,0)$.
\end{lemma}

\begin{remark}
  The notation $H^*(j)$ means the pull-back of $j$ seen as a function
  on $\overline{\mathcal{M}}^1$, by the algebraic morphism $H:
  \overline{\mathcal{M}}^1 \rightarrow 
  \overline{\mathcal{M}}^1$. We prefer it instead of $H(j)$, in order
  not to get 
  confused in the next section by a notation of the type $H(h)$, which
  is not simply obtained from $H(j)$ by replacing $j$ with $h$. 
\end{remark}

In order to find the unknown coefficients $\alpha, \beta, \gamma,
\delta $, we look first for the order of the zero $j= 2^8\cdot 3^3$ and
of the pole $j= 0$ of $H(j)$. Writing $M:=m^3$, relation
(\ref{commdiag}) becomes :
  \begin{equation} \label{eqM}
      H\bigg(27\cdot \dfrac{M(M+8)^3}{(M-1)^3}\bigg)=
        \dfrac{(4-M)^3}{M^2}\cdot 
          \bigg(\dfrac{(M-4)^3 -6^3 M^2}{(M-4)^3 +3^3 M^2}\bigg)^3.
  \end{equation}

By factoring $H^*(j) =\dfrac{(j-2^8\cdot 3^3)^k}{j^l}\cdot
K(j)$, with $K(2^8\cdot 3^3) \in \mathbb{C}\setminus
\{ 0, \infty\}$ and $K(0) \in \mathbb{C}\setminus
\{ 0, \infty\}$, we deduce easily that $k=3$ and $l=2$. Combining this
with Lemma \ref{prelform}, we see that:
 $$H^*(j) =\dfrac{1}{\gamma}\cdot \dfrac{(j - 2^8\cdot 3^3)^3}{j^2}.$$
By equating  the dominating coefficients of both sides of (\ref{eqM}),
when $M\rightarrow \infty$, we deduce that $\gamma=-27$. 

We have got like this the desired expression of the hessian dynamical
system:

\begin{theorem}  \label{exphess}
   One has the following expression of the hessian dynamical system in
   terms of the parameter $j$, on the compactified modular curve
   $\overline{\mathcal{M}}_1$:
$$H^*(j)= -\frac{1}{3^3}\cdot \frac{(j - 2^8\cdot 3^3)^3}{j^2}.$$
\end{theorem}

\medskip
\section{The associated sequence of dessins d'enfants}
  \label{Secdes} 

Let us make the following change of variable on the modular curve
$\mathcal{M}_1$:
\begin{equation} \label{chmod}
   j= 2^6\cdot 3^3 \cdot h
\end{equation}

From Proposition \ref{exphess}, we deduce the following expression of
the hessian dynamical system in terms of the variable $h$:
\begin{equation} \label{eprh}
   H^*(h)= -\dfrac{1}{27}\cdot \dfrac{(h-4)^3}{h^2}
\end{equation}
which shows that:
 $$\dfrac{dH^*(h)}{dh}=-\dfrac{1}{27}\cdot
 \dfrac{(h-4)^2(h+8)}{h^3}.$$

We deduce from this immediately:

\begin{proposition} \label{locim}
   Set-theoretically, the critical locus of
   $H:\mathbb{P}^1_{[1:h]}\rightarrow \mathbb{P}^1_{[1:h]}$ is equal to 
   $$\mathrm{Crit}(H)= \{h=4, h=-8, h=0\}.$$ 
   The critical image of $H$, also called its discriminant locus,  is equal to:
   $$\mathrm{\Delta}(H)= \{h=0, h=1, h =\infty\}.$$
   Seen as divisors, the critical fibers of $H$ are:
   $$\begin{array}{l}
        \mathrm{div}(H^*(h))=3(h=4),\\
        \mathrm{div}(H^*(h-1))=2(h=-8) + (h=1),\\
        \mathrm{div}(H^*(1/h))= 2(h=0)+ (h=\infty), 
     \end{array}$$
   where $(h=a)$ denotes the point of $\mathbb{P}^1_{[1:h]}$ where the
   rational function $h$ takes the value $a \in \mathbb{C} \cup
   \{\infty\}$. 
\end{proposition}

The previous proposition explains why we have chosen the change of
variable (\ref{chmod}): in order to get as critical image the set
$\{0,1,\infty\}$ of values of the working parameter.

The elliptic curves corresponding to the critical values $h=0$ and
$h=1$ of $H$ inside $\mathcal{M}_1=\mathbb{C}_h$ are exactly those
with special symmetry, as shown by the following proposition (see
Hartshorne \cite[page 321]{H 77}): 

\begin{proposition} \label{speciell}
  Let $E$ be an elliptic curve over $\mathbb{C}$. Denote by $G_E$ the
  group of automorphisms of $E$ leaving a base point fixed. Then $G_E$ is a
  finite group of order:

$\bullet$ $2$ if $j(E) \notin \{ 0, 1728\} \Leftrightarrow h
\notin\{0, 1\}$.

$\bullet$ $4$ if $j(E)=1728 \Leftrightarrow h=1$.

$\bullet$ $6$ if $j(E)=0 \Leftrightarrow h=0$.

\end{proposition}

By Proposition \ref{locim}, we see that the cardinal of the
discriminant set $\Delta(H)$ 
of the hessian dynamical system is equal to $3$ when we look at $H$  as
a ramified covering of $\mathbb{P}^1$. By the
work \cite{B 80} of Belyi, we know that ramified covers of
$\mathbb{P}^1$ with $3$ critical values are particularly important
from the arithmetical viewpoint (see also Zapponi \cite{Z
  03}). Following this last reference, let us recall the notion of
\emph{dessin d'enfant} associated to such a map, introduced initially
at the suggestion of  Grothendieck \cite[section 3]{G 97}.

Let $\psi: C \rightarrow \mathbb{P}^1:= \mathbb{C}\cup \{ \infty\}$ be
a holomorphic map from a compact Riemann surface $C$ to
$\mathbb{P}^1$. Denote by $\Gamma$ the preimage
$\psi^{-1}([0,1])$. Color the vertices of $\Gamma\cap \psi^{-1}(0)$ in
\emph{black} and those of   $\Gamma\cap \psi^{-1}(1)$ in
\emph{white}. Moreover, order cyclically the germs of edges starting
from each vertex of $\Gamma$ as they occur when one turns positively
with respect to the canonical orientation defined by the complex
structure of $C$. 

\begin{definition} \label{dessin}
   The graph $\Gamma$ with colored vertices and cyclically ordered
   germs of edges obtained as explained before is called the
   \textbf{dessin d'enfant} associated to $\psi$.
\end{definition}

The point of this definition is that this dessin (a purely
\emph{topological} object) encodes completely up to isomorphisms the
map $\psi$ (a \emph{holomorphic} object).

\begin{remark}
    If one has a map $\psi:C\rightarrow P$ where $P$ is isomorphic to
    $\mathbb{P}^1$ and the discriminant set has cardinal equal to $3$,
    one has to choose which points between the $3$ critical values are
    to be identified with $0$ and $1$ in order to define the dessin
    d'enfant associated to $\psi$. In our case, the point $\infty$ is
    determined geometrically as the point at infinity of the moduli
    space  $\mathcal{M}_1$. 
\end{remark}

Proposition \ref{locim} shows that $H(\Delta(H))\subset
\Delta(H)$. More precisely, $0 \rightarrow \infty, 1 \rightarrow 1,
\infty \rightarrow \infty$. This implies: 

\begin{proposition}\label{iterates}
 All the iterates $H^{(n)}:= \underbrace{H \circ \cdots \circ H}_{n\:
   \mbox{times}}$ 
 of $H$, where
$n \geq 1$, are also ramified covers of $\overline{\mathcal{M}}_1:=
\mathbb{P}^1_{[1:h]}$, with ramification set $\{0,1,\infty\}$. 
\end{proposition}

We deduce that each iterate has an associated dessin d'enfant
$\Gamma_n$. As the map $H$ is defined over $\mathbb{Q}$, all its
iterates have the same property, which shows that all the dessins 
$\Gamma_n$ are fixed under the natural action of the absolute Galois 
group $\mbox{Gal}(\overline{\mathbb{Q}}/ \mathbb{Q})$.

\medskip

\section{An algorithm for constructing the sequence 
  of dessins d'enfants}  \label{Secomb}

We want now to understand  how evolves the sequence of dessins
d'enfants 
$(\Gamma_n)_{n \geq   1}$. The graph $\Gamma_n$ embeds into the
preimage $(H^{(n)})^{-1}(\mathbb{R}_h \cup \{\infty\})$ of the real
projective line of the parameter $h$. 

From the holomorphic view-point, the real projective line 
  $\mathbb{R}_h \cup \{\infty\}$ is
canonically determined by the dynamical system, as the unique circle
$C_{\Delta}$ 
contained in the smooth projective rational curve
$\overline{\mathcal{M}}_1= \mathbb{P}_{[1:h]}^1$, and which contains the
discriminant set $\Delta(H)$.  

Denote:
  $$G_n := (H^{(n)})^{-1}(C_{\Delta}).$$
From the topological view-point, $G_n$ is a graph embedded into
$\overline{\mathcal{M}}_1$, its vertices being the preimages of $\{0,1,\infty\}
\in \mathbb{P}_{[1:h]}^1= \overline{\mathcal{M}}_1$. We decorate the
edges of $G_n$ in three different
ways, according to the real interval $(\infty,0), (0,1), (1, \infty)$
which is their image under $H$. We orient them with the lift by $H$ of
the natural orientation of $\mathbb{R}_h$ from negative to positive
numbers. The drawing convention we choose is indicated in Figure 1.

 {\tt    \setlength{\unitlength}{0.92pt}}
\begin{figure} \label{Convaxis}
   \epsfig{file=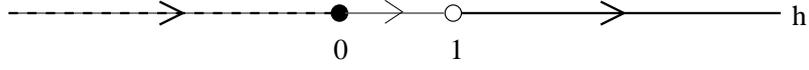, height= 8 mm}
   \caption{The drawing convention}
\end{figure}

{\tt    \setlength{\unitlength}{0.92pt}}
\begin{figure} \label{Preimreel1}
   \epsfig{file=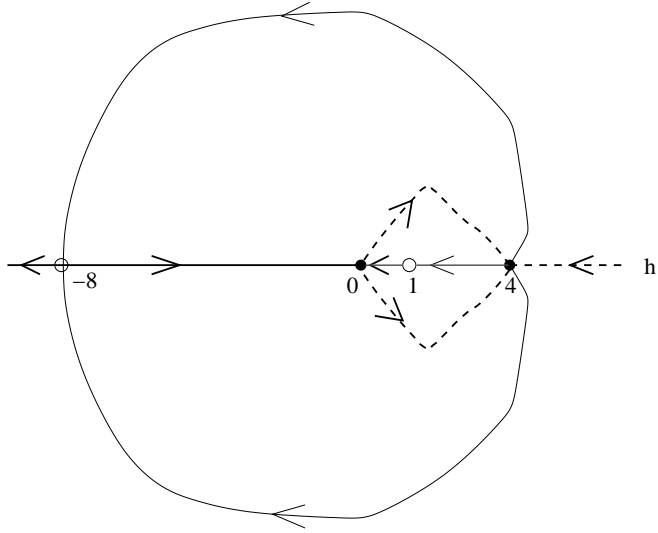, height= 70 mm}
   \caption{The preimage by $H$ of the real line}
\end{figure}

{\tt    \setlength{\unitlength}{0.92pt}}
\begin{figure} \label{Dessin1}
   \epsfig{file=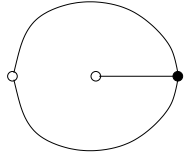, height= 20 mm}
   \caption{The dessin d'enfant $\Gamma_1$ associated to $H$}
\end{figure}

Write $$x:=\mathrm{Re}(h),\: y:= \mathrm{Im}(h).$$ Then the equation 
$\mathrm{Im}(H^*(h))=0$ of $H^{-1}(C_{\Delta})$ becomes:
$$\mathrm{Im}\bigg[ -\dfrac{1}{27}\cdot \dfrac{(x+iy-4)^3}{(x+iy)^2}\bigg]=0.$$
After a few computations we get the equation:
\begin{equation} \label{eqpreim}
   y \cdot\bigg[((x-4)^2 +y^2)^2+16(x-4)((x-4)^2 +y^2)+
   16(3(x-4)^2-y^2)\bigg]=0.
\end{equation}

This shows that $G_1= H^{-1}(C_{\Delta})$ is the union of the
real axis of the variable $h$ and a singular quartic curve, which has
as only singularity a real node at $h=4$. The union of the two curves
is drawn in Figure 2. 

By looking inside Figure 2 at $H^{-1}([0,1])$, we deduce
that the dessin d'enfant $\Gamma_1$ 
associated to $H$ is as indicated in Figure 3.

\medskip

Consider the map represented in Figure 4. 
In it, $T$ is a compact affine triangle with the vertices denoted
$0,1,\infty$ and $H_{PL}:T \rightarrow T$ is a continuous
piecewise-linear map, which is an affine homeomorphism onto $T$ in
restriction to 
each of the three closed triangles into which $T$ is triangulated. Their
vertices are midpoints of the edges of $T$.  The distinct edges
of the $1$-skeleton at the source are decorated as their images by
$H_{PL}$. 

{\tt    \setlength{\unitlength}{0.92pt}}
\begin{figure} \label{Preimage1}
   \epsfig{file=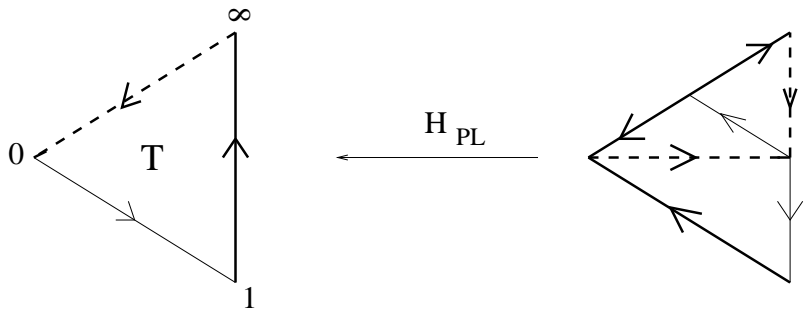, height= 40 mm}
   \caption{The combinatorics of the triangulation $\mathcal{T}_1$}
\end{figure}

{\tt    \setlength{\unitlength}{0.92pt}}
\begin{figure} \label{Preimage2}
   \epsfig{file=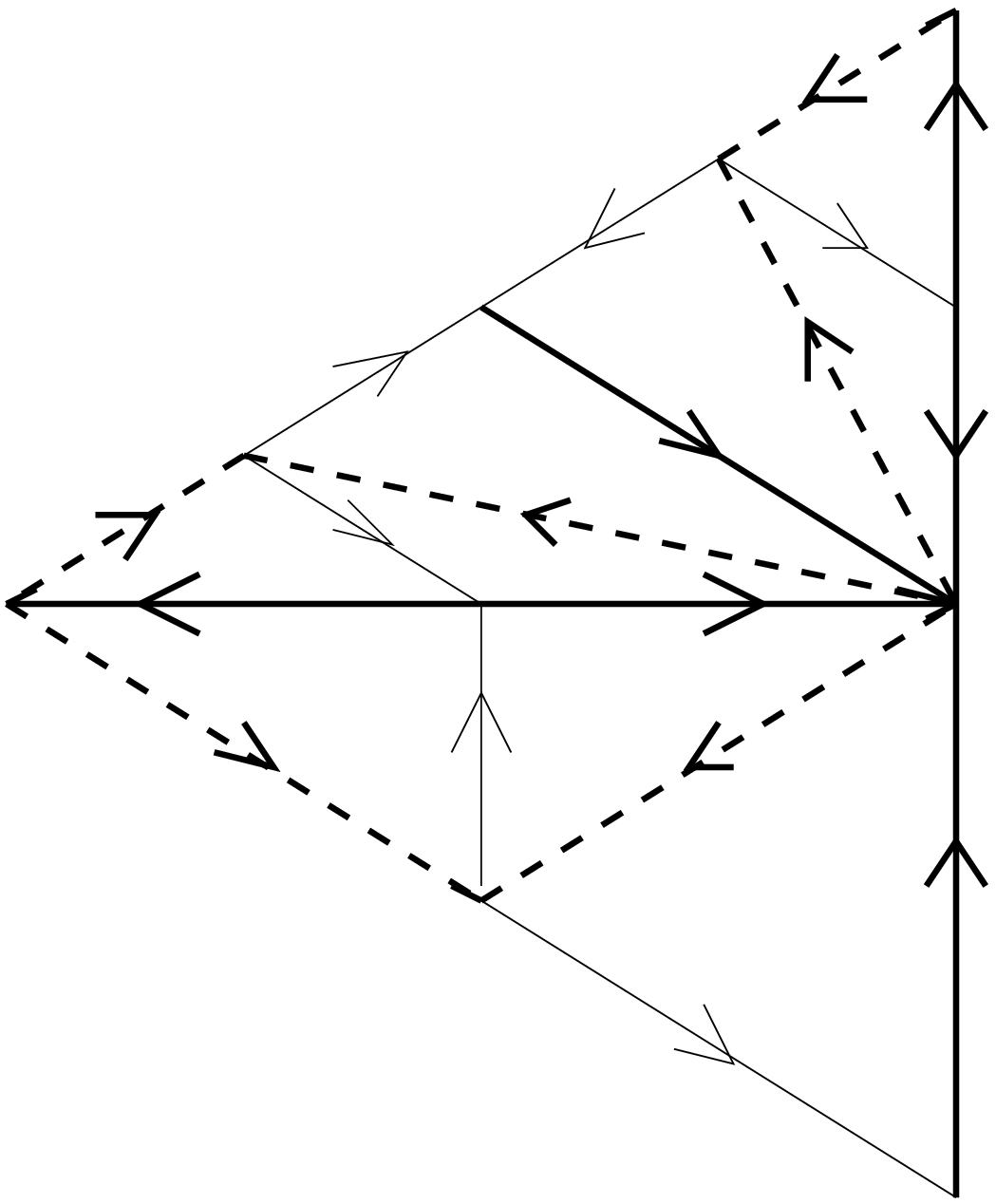, height= 50 mm}
   \caption{The combinatorics of the triangulation $\mathcal{T}_2$}
\end{figure}

{\tt    \setlength{\unitlength}{0.92pt}}
\begin{figure} \label{Preimage3}
   \epsfig{file=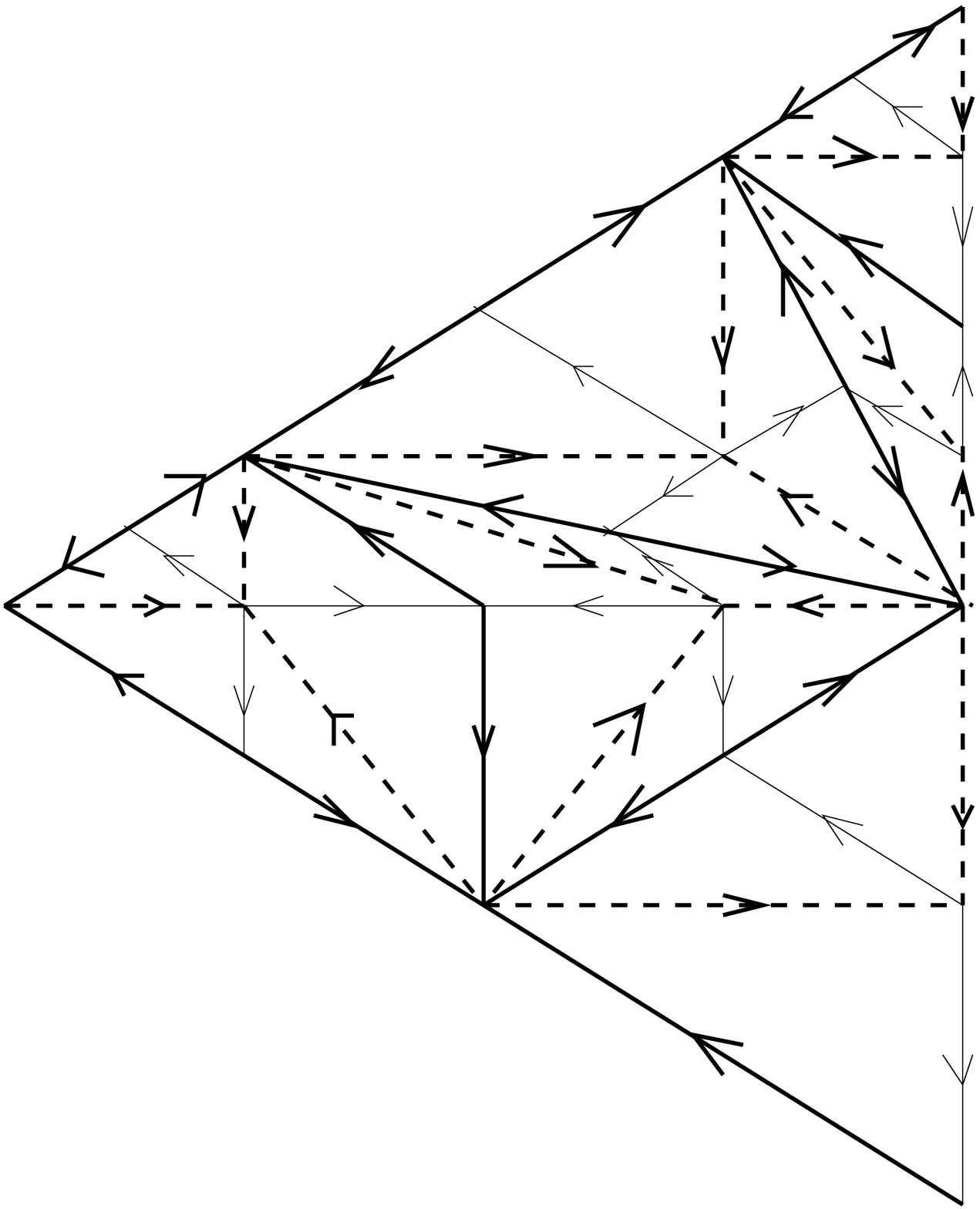, height= 80 mm}
   \caption{The combinatorics of the triangulation $\mathcal{T}_3$}
\end{figure}

{\tt    \setlength{\unitlength}{0.92pt}}
\begin{figure} \label{Preimage4}
   \epsfig{file=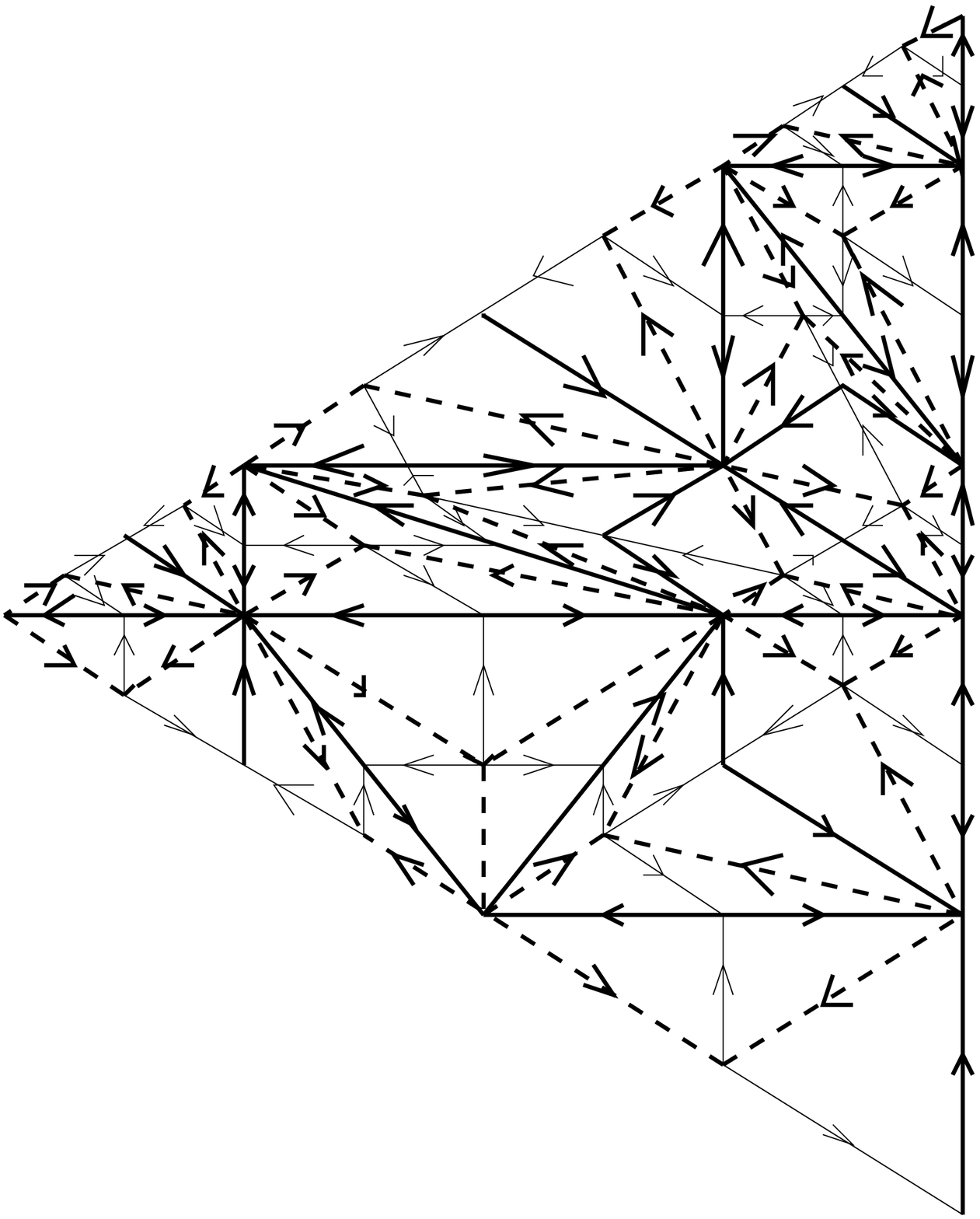, height= 110 mm}
   \caption{The combinatorics of the triangulation $\mathcal{T}_4$}
\end{figure}

{\tt    \setlength{\unitlength}{0.92pt}}
\begin{figure} \label{Dessindenf2}
   \epsfig{file=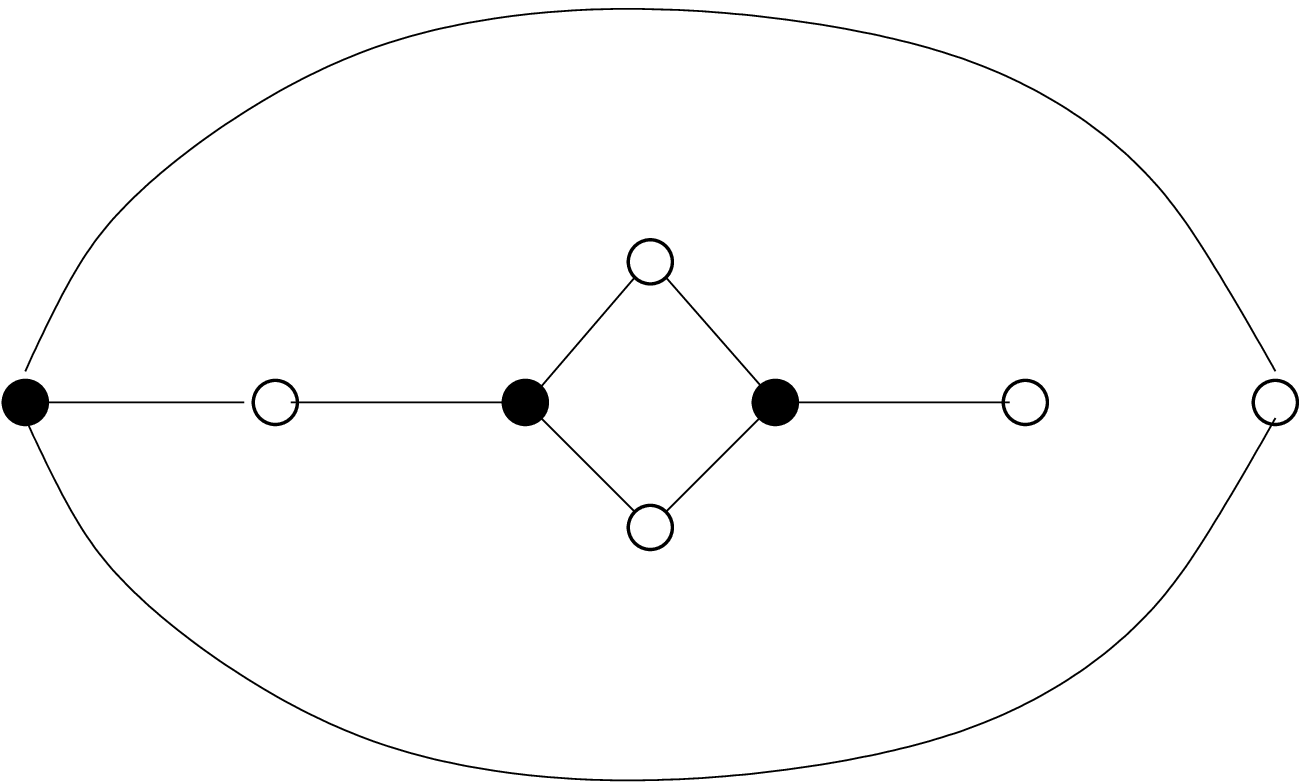, height= 30 mm}
   \caption{The dessin d'enfant $\Gamma_2$ associated to $H^{(2)}$}
\end{figure}

{\tt    \setlength{\unitlength}{0.92pt}}
\begin{figure} \label{Dessindenf3}
   \epsfig{file=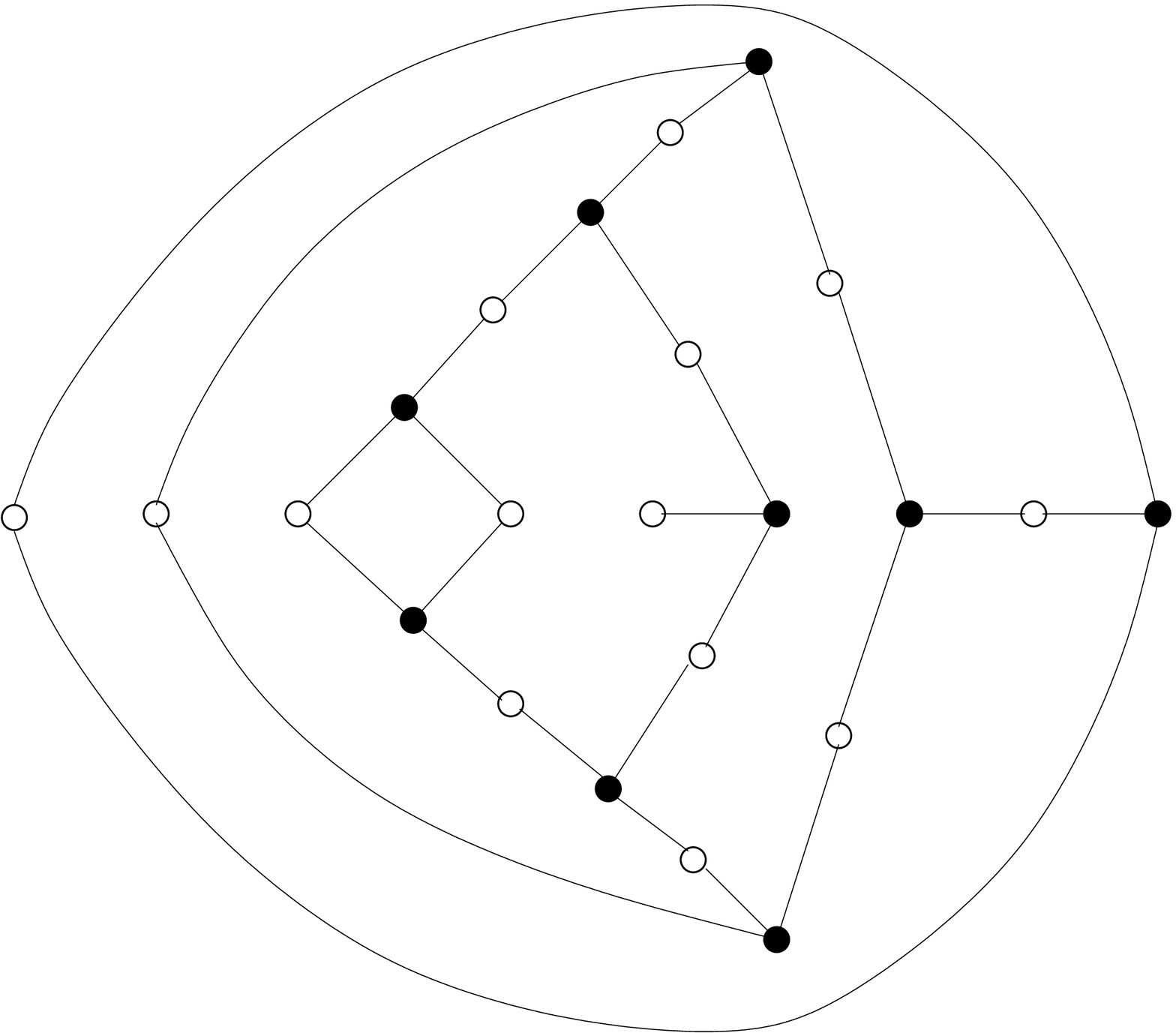, height= 50 mm}
   \caption{The dessin d'enfant $\Gamma_3$ associated to $H^{(3)}$}
\end{figure}

{\tt    \setlength{\unitlength}{0.92pt}}
\begin{figure} \label{Dessindenf4}
   \epsfig{file=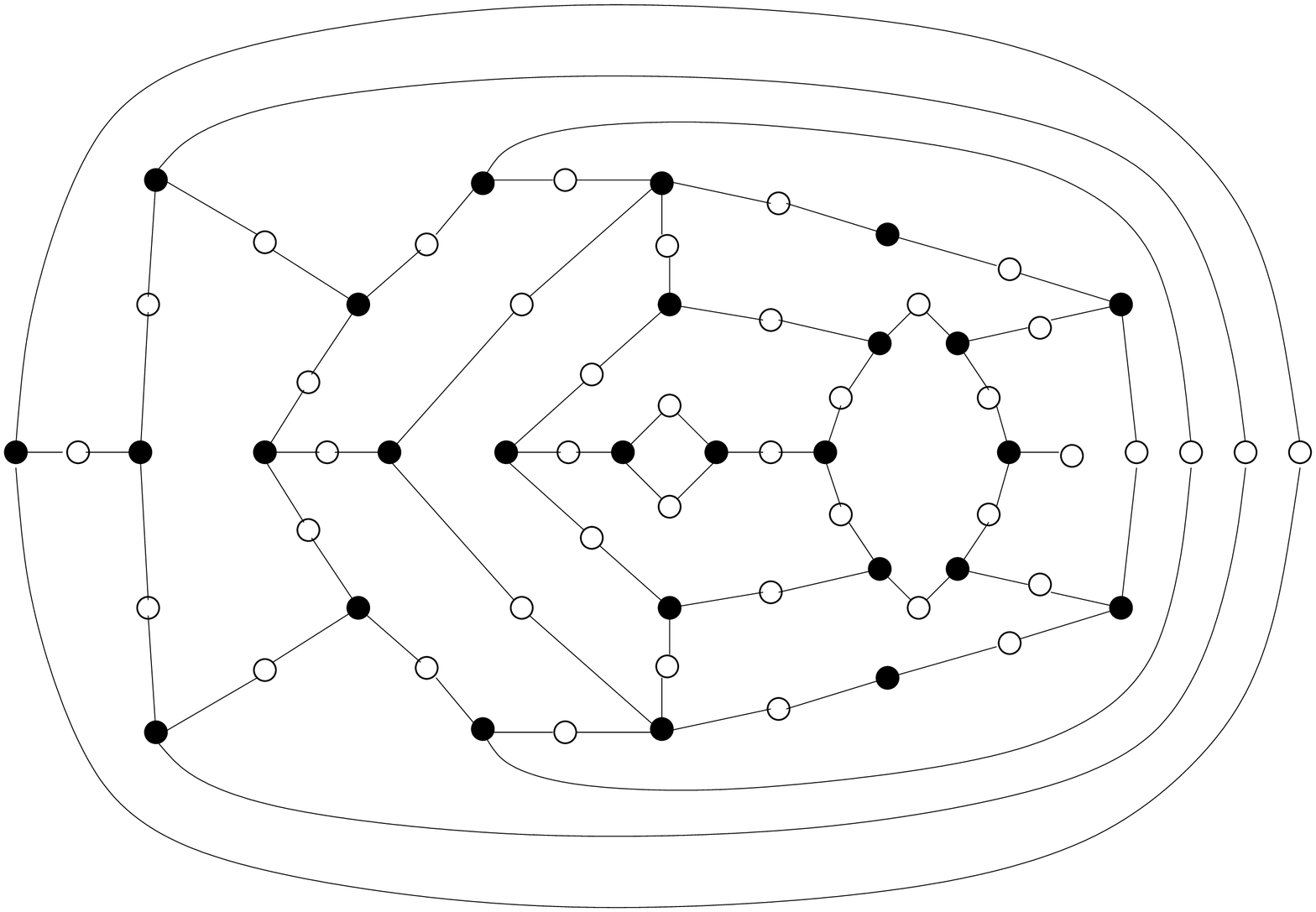, height= 70 mm}
   \caption{The dessin d'enfant $\Gamma_4$ associated to $H^{(4)}$}
\end{figure}

We see that the topology of the map is completely described by the
triangulation of the source triangle and by the decorations of the
edges. In the same way, we can describe the topology of $H_{PL}^{(n)}, \:
\forall \: n\geq 1$, by a decorated triangulation:
$$\mathcal{T}_n:=(H_{PL}^{(n)})^{-1}(\partial T).$$ 
The convention for the decoration of the edges of $\mathcal{T}_n$ is
the same as the one used for $G_n$. 
We have the following algorithm for the construction of the sequence $(
\mathcal{T}_n)_{n \geq 1}$:
\medskip

$\bullet$ Start with $H_{PL}$ and call the triangulation at the source
  $\mathcal{T}_1$. 

$\bullet$  Given the triangulation $\mathcal{T}_n$, construct
   $\mathcal{T}_1$  inside each triangle
   of $\mathcal{T}_n$. Call the triangulation obtained like this
   $\mathcal{T}_{n+1}$.

Let $\Phi_n$ be the subgraph of the $1$-skeleton $\mathcal{T}_n^1$ of the
  triangulation $\mathcal{T}_n$ obtained as the union of all edges
  whose image by $H_{PL}^{(n)}$ is the edge $[0,1]$ of $T$. Decorate the
  preimages of $0$ in black and those of $1$ in white. Then take the
  double $(T, \mathcal{T}_n^1, \Phi_n) \cup_{\partial T} (T,
  \mathcal{T}_n^1, \Phi_n)$ of the triple $(T,\mathcal{T}_n^1, 
  \Phi_n)$ of topological spaces, that is glue two copies of it by
  the identity on 
  $\partial T$. Orient arbitrarily the sphere $T \cup_{\partial T} T$
  and orient cyclically the germs of edges of $\Phi_n\cup_{\partial T}
  \Phi_n$ around each vertex according to this orientation. Denote
  $\Sigma:= T \cup_{\partial T} T, \: \tilde{\mathcal{T}}_n^1:= 
  \mathcal{T}_n^1\cup _{\partial T \:\cap \:
    \mathcal{T}_n^1}\mathcal{T}_n^1, \: \tilde{\Phi}_n:= \Phi_n
  \cup_{\partial T \: \cap \: \Phi_n} \Phi_n$. 

\begin{proposition} \label{embgraph}
  The decorated  graph $\tilde{\Phi}_n$ is homeomorphic (respecting
  the decorations) to the dessin  d'enfant $\Gamma_n$ of  $H^{(n)}$. 
\end{proposition}

\medskip

In Figures 5,6 and 7 we have drawn the decorated triangulations
$\mathcal{T}_2, \mathcal{T}_3$ and $\mathcal{T}_4$. By Proposition
\ref{embgraph}, one can extract
easily from them the dessins d'enfants $\Gamma_2, \Gamma_3$ and
$\Gamma_4$, illustrated in Figures 8,9 and 10 respectively. We have
chosen a homeomorphism between $T$ and a half-plane, sending the
vertex $\infty$ of $T$ to infinity, and then we have drawn the double
$\tilde{\Phi}_n$ of $\Phi_n$ just by gluing the reflexion of $\Phi_n$
with respect to the border-line of the half-plane.

{\tt    \setlength{\unitlength}{0.92pt}}
\begin{figure} \label{Eucl1}
   \epsfig{file=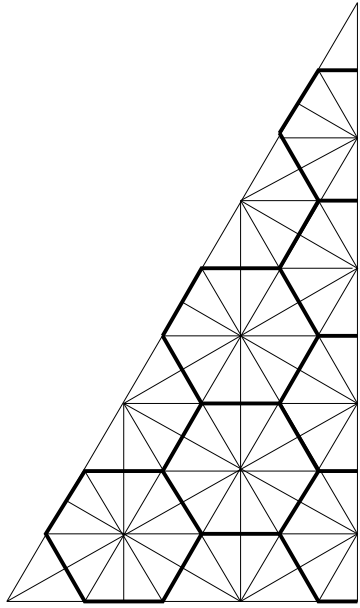, height= 80mm}
   \caption{Another picture of $\mathcal{T}_4$}
\end{figure}

\begin{remark} \label{add}
(Added in proofs) After my exposition of the results of
  this paper in the University of Toulouse in December 2008, X. Buff
  and A. Ch{\'e}ritat 
  proved easily using general theory that the map $H$ is a  so-called
  \emph{Latt{\`e}s map} (see \cite{L 18} and \cite{M 04}). More
  precisely, it is conjugate to the map induced on $\mathbb{C}/G$ by the
  multiplication $\cdot i \sqrt{3}: \mathbb{C}\rightarrow
  \mathbb{C}$. Here $G$ denotes the group of automorphisms of
  $\mathbb{C}$ generated by the translations $z \rightarrow z+1, \: z
  \rightarrow z + \epsilon$ and the rotation $z \rightarrow \epsilon z$, where
  $\epsilon := e^{\frac{i \pi}{3}}$. This leads to a nice
  description of the sequences $(\mathcal{T}_n)_n$ and $(\Gamma_n)_n$
using euclidean geometry, in the spirit of Thurston's philosophy that
topology is best understood using adapted geometric structures. The
point here is that one should consider on $\mathbb{C}/G$ the structure
of euclidean orbifold, quotient of the euclidean structure of $\mathbb{C}$ by
$G$, which is a group of isometries for it. In this geometry, the
figures of this paper appear in a way which is completely natural
within the framework of euclidean geometry.  Namely, the iterative
operation drawn in Figure 4 should be replaced by one in which the
triangle $T$ is similar to each one of the 3 subtriangles created by
the subdivision, the similarities respecting the decorations of the
vertices by the symbols $0, 1, \infty$. This forces the angles of $T$
to be of $60^{\circ}, 90^{\circ}$ and $30^{\circ}$ respectively. In
Figure 11 is drawn in this way $\mathcal{T}_4$, as well as the
union of those edges allowing to construct $\Gamma_4$ by
doubling. Compare with Figure 7 !
\end{remark}

\medskip


\medskip

{\small


\begin{thebibliography}{99}

\bibitem{AD 06} Artebani, M., Dolgachev, I. \textit{The Hesse pencil
    of plane cubic curves.} arXiv:math/0611\\590. To appear in
  l'Enseignement Math{\'e}matique. 

\bibitem{B 80} Belyi, G. V., \textit{On Galois extensions of a maximal
    cyclotomic field.} Mathematical USSR Izvestija \textbf{14} (1980),
  247-256. 

\bibitem{BK 86}Brieskorn, E., Kn{\"o}rrer, H. \textit{Plane algebraic
  curves.} Translated from the German by John Stillwell. Birkh{\"a}user
  Verlag, Basel, 1986.  

\bibitem{E 03} Edidin, D. \textit{What is... a Stack?} Notices of the
  AMS \textbf{50}, no. 4 (2003), 458-459.

\bibitem{G 97} Grothendieck, A. \textit{Sketch of a programme.} In
  \textit{Geometric Galois Actions.} Schneps, L., Lochak, P. eds.,
  London Math. Soc. Lecture Note Ser. \textbf{242}, Cambridge
  Univ. Press, 1997. English translation of \textit{Esquisse d'un
    programme.} Preprint, 1984. 

\bibitem{H 77} Hartshorne, R. \textit{Algebraic Geometry.} Springer,
  1977. 

\bibitem{H 26} Hollcroft, T. \textit{Harmonic cubics.}
  Ann. Math. \textbf{27} (1926), 568-576. 

\bibitem{L 18} Latt{\`e}s, S. \textit{Sur l'it{\'e}ration des
    substitutions rationnelles et les fonctions de Poincar{\'e}.}
  C. R. Acad. Sci. Paris \textbf{166} (1918), 26-28. 

\bibitem{M 04} Milnor, J. \textit{On Latt{\`e}s maps.} In
  \textit{Dynamics on the Riemann sphere.} 9-43. Eur. Math. Soc.,
  Z{\"u}rich, 2006.  

\bibitem{P 00} Pilgrim, K.M. \textit{Dessins d'enfants and Hubbard
    trees.} Ann. Sci. Ec. Norm. Sup. \textbf{33} (2000), 671-693. 

\bibitem{Z 03} Zapponi, L. \textit{What is... a Dessin d'Enfant?}
  Notices of the  AMS \textbf{50}, no. 7 (2003), 788-789.
\end{thebibliography} }
\end{document}